\documentclass[12pt]{amsart}
\usepackage{amssymb}
\usepackage{amsmath}
\usepackage{amsfonts}
\usepackage{graphicx}
\usepackage{color}
\usepackage[onehalfspacing]{setspace}
\usepackage{hyperref}

\numberwithin{equation}{section}
\theoremstyle{plain}

\input{tcilatex}

\begin{document}
\title[Minimax Estimation and Discontinuity]{Minimax Estimation of
Nonregular Parameters and Discontinuity in Minimax Risk}
\date{March 10, 2014.}
\author[K. Song]{Kyungchul Song}
\address{Vancouver School of Economics, University of British Columbia, 997
- 1873 East Mall, Vancouver, BC, V6T 1Z1, Canada}
\email{kysong@mail.ubc.ca}

\begin{abstract}
{\footnotesize When a parameter of interest is nondifferentiable in the
probability, the existing theory of semiparametric efficient estimation is
not applicable, as it does not have an influence function. Song (2014)
recently developed a local asymptotic minimax estimation theory for a
parameter that is a nondifferentiable transform of a regular parameter,
where the nondifferentiable transform is a composite map of a continuous
piecewise linear map with a single kink point and a translation-scale
equivariant map. The contribution of this paper is two fold. First, this
paper extends the local asymptotic minimax theory to nondifferentiable
transforms that are a composite map of a Lipschitz continuous map having a
finite set of nondifferentiability points and a translation-scale
equivariant map. Second, this paper investigates the discontinuity of the
local asymptotic minimax risk in the true probability and shows that the
proposed estimator remains to be optimal even when the risk is locally
robustified not only over the scores at the true probability, but also over
the true probability itself. However, the local robustification does not
resolve the issue of discontinuity in the local asymptotic minimax risk.}

{\footnotesize \ }

{\footnotesize \noindent \textsc{Key words.} Nonregular Parameters;
Semiparametric Efficiency; Local Asymptotic Minimax Estimation;
Translation-Scale Equivariant Maps \newline
}

\end{abstract}

\maketitle

\section{Introduction}

Statistical inference on a parameter begins by choosing an appropriate
estimator. For a finite dimensional parameter defined under local asymptotic
normal experiments, it has become nearly a standard practice in statistics
and econometrics to establish the optimality of an estimator through
semiparametric efficiency, where optimality is expressed as a variance
bound, and an estimator is taken to be optimal if it is asymptotically
normal with its asymptotic variance achieving the bound. The literature
along this approach is vast in statistics and econometrics.

A mathematical analysis of asymptotic optimal inference began with the
famous paper by Wald (1943). While the approach of local asymptotic minimax
estimation has appeared in the previous literature (e.g. Le Cam (1953) and
Chernoff (1956)), major breakthroughs were made by H\'{a}jek (1972) and Le
Cam (1972). Koshevnik and Levit (1976), Pfanzagl and Wefelmeyer (1982),
Begun, Hall, Huang and Wellner (1983) and Chamberlain (1986) extended
asymptotic efficient estimation to nonparametric and semiparametric models.
See also van der Vaart (1988, 1991) for further developments in this
direction, and Newey (1990) for results that are relevant to econometrics. A
general account of this approach is found in monographs such as Bickel,
Klaassen, Ritov, and Wellner (1993), and in later chapters of van der Vaart
and Wellner (1996) and van der Vaart (1998). While the references so far
mostly focus on local asymptotic normal experiments (as this paper does),
optimal estimation theory in local asymptotic mixed normal experiments has
also received attention in the literature. See Jeganathan (1982) and Basawa
and Scott (1983). In econometrics, Phillips (1991) developed optimal
inference theory for cointegrating regression models using the framework of
local asymptotic mixed normal experiments.

The existing notion of semiparametric efficiency is not directly applicable,
when the parameter is not differentiable in the probability that identifies
the parameter. Nondifferentiable parameters do not merely constitute a
pathological case, for one can easily encounter such a parameter when the
parameter (denoted by $\theta \in \mathbf{R}$) is defined through a
nondifferentiable transform of another parameter vector, say, $\mathbf{\beta 
}\in \mathbf{R}^{d}$. For example, the parameter of interest might take the
form of $\theta =\max \{\beta _{1},\cdot \cdot \cdot ,\beta _{d}\}$ or $%
\theta =\min \{\beta _{1},\cdot \cdot \cdot ,\beta _{d}\}$, where $\beta
_{1},\cdot \cdot \cdot ,\beta _{d}$ are average treatment effects from
different treatment regimes or mean squared prediction errors from different
predictive models, or boundaries of multiple intervals. (As for the last
example, Chernozhukov, Lee and Rosen (2013) called the parameter an
intersection bound. Examples of such bounds are found in Haile and Tamer
(2003) and Manski and Tamer (2002) among many others.) The difficulty with
estimation theory for such nondifferentiable parameters is emphasized by
Doss and Sethuraman (1988). See Hirano and Porter (2010) for a general
impossibility result for such parameters.

This paper focuses on the problem of optimal estimation when the parameter
is nondifferentiable in a particular way. More specifically, this paper
focuses on a parameter of interest in $\mathbf{R}$ which takes the following
form:%
\begin{equation}
\theta =(f\circ g)(\mathbf{\beta }),  \label{param2}
\end{equation}%
where $\mathbf{\beta }\in \mathbf{R}^{d}$ is a regular parameter for which a
semiparametric efficiency bound is well defined, $g$ is a translation-scale
equivariant map, and $f$ is a continuous map that is potentially
nondifferentiable. While the paper focuses on this particular way that
nondifferentiability arises, it accommodates various nondifferentiable
parameters that are relevant in empirical researches (See Song (2014) for
examples.)

A recent work by the author (Song (2014)) considers the case of $f$ being a
continuous piecewise linear map \textit{with a single kink point}, and has
demonstrated that the existing semiparametric efficient estimation can be
extended to this case of nonregular parameter $\theta $ through a local
asymptotic minimax approach. While the result applies to various examples of
nonregular parameters used in econometrics, the restriction on $f$ excludes
some interesting examples. For example, one might be interested in an
optimal policy parameter that is censored on both upper and lower bounds,
say, due to constraints on resources or in implementation.

This paper generalizes the theory to the case where $f$ is Lipschitz
continuous yet potentially nondifferentiable at a finite number of points.
Similarly, as in Song (2014), it turns out that the local asymptotic minimax
estimator takes the following form:%
\begin{equation}
\hat{\theta}_{mx}\equiv f\left( g(\mathbf{\tilde{\beta}})+\frac{\hat{c}%
^{\ast }}{\sqrt{n}}\right) ,  \label{oe}
\end{equation}%
where $\hat{c}^{\ast }$ is an \textit{optimal bias adjustment term}, and $%
\mathbf{\tilde{\beta}}$ is a semiparametrically efficient estimator of $%
\mathbf{\beta }$. The optimal bias adjustment term can be determined by
simulating the local asymptotic minimax risk.

Some researches in the literature have suggested various methods of bias
adjustment and reported improved performances. (See for example Haile and
Tamer (2003), and Chernozhukov, Lee, and Rosen (2013).) The approach of Song
(2014) and this paper is distinct in the sense that it determines the
optimal bias adjustment explicitly through theory of local asymptotic
minimax estimation.

The resulting local asymptotic minimax risk for the kind of nonregular
parameters considered in this paper is discontinuous in the underlying true
probability in general. To appreciate the meaning of this discontinuity, it
is worth recalling that the classical local asymptotic minimax risk approach
imposes local uniformity over parametric submodels passing a fixed true
probability. This local uniformity eliminates superefficient estimators such
as Hodges estimator which is known to exhibit poor finite sample
performance. (See Le Cam (1953) for a formal treatment of Hodges
superefficient estimator. See also Weiss and Wolfowitz (1966)).) In
classical estimation theory for regular parameters, the local asymptotic
minimax risk is continuous in the true probability. This feature stands in
contrast with the local asymptotic minimax risk in this paper which is
discontinuous in the true probability.

When the asymptotic distribution of a test statistic or an estimator
exhibits discontinuity in the underlying true probability, it is a common
practice to consider an alternative asymptotic theory along a sequence of
probabilities local around the true probability. Mostly, this alternative
asymptotics involves a localization parameter which continuously "bridges"
two distributions across the discontinuity point. The most common example of
this approach is local asymptotic power analysis in hypothesis tests, where
one adopts a sequence of probabilities that converge to a probability that
belongs to the null hypothesis. A similar approach is found in local to
unity models (Stock (1991)), weakly identified models (Staiger and Stock
(1997)), and more recently, models of various moment inequality restrictions
(Andrews and Guggenberger (2009)) among many others.

To deal with this issue of discontinuity in the local asymptotic minimax
risk, this paper introduces a local robustification of the risk, where the
risk is further robustified against a local perturbation of the true
probability. Somewhat unexpectedly, the local asymptotic minimax risk
remains unchanged after this local robustification of the risk. On the one
hand, this means that the local asymptotic minimax estimator in (\ref{oe})
retains its optimality under this robustification, as long as the efficient
estimator $\mathbf{\tilde{\beta}}$, after location-scale normalization,
converges in distribution uniformly over the true probabilities. On the
other hand, the discontinuity of the risk in the true probability is not
resolved by the local robustification approach. Hence there may be a gap
between the finite sample risk and its asymptotic version and the gap does
not close uniformly over all the probabilities even in the limit. It remains
an open question whether this renders the whole apparatus of the local
asymptotic minimax theory dubious in practice, when the parameter is
nondifferentiable.

The rest of the paper is structured as follows. In Section 2, the paper
defines the scope of this paper by introducing assumptions about $f$, $g$, $%
\mathbf{\beta }$, and the set of underlying probabilities that identify $%
\mathbf{\beta }$. In Section 3, the paper gives a characterization of local
asymptotic minimax risk, and proposes a general method to construct a local
asymptotic minimax estimator. In Section 4, the paper considers a local
robustification of the local asymptotic minimax risk, and shows that the
results of Section 3 mostly remain unchanged. Section 5 concludes the paper.
The mathematical proofs of the paper's results appear in the Appendix.

A word of notation. Let $\mathbf{1}_{d}$ be a $d\times 1$ vector of ones
with $d\geq 2$. For a vector $\mathbf{x}\in \mathbf{R}^{d}$ and a scalar $c$%
, we simply write $\mathbf{x}+c=\mathbf{x}+c\mathbf{1}_{d}$, or write $%
\mathbf{x}=c$ instead of $\mathbf{x}=c\mathbf{1}_{d}$. For $\mathbf{x}\in 
\mathbf{R}^{d}$, the notation $\max (\mathbf{x})$ (or $\min (\mathbf{x})$)
means the maximum (or the minimum) over the entries of the vector $\mathbf{x}
$. We let $\mathbf{\bar{R}}=[-\infty ,\infty ]$ and view it as a two-point
compactification of $\mathbf{R}$, and let $\mathbf{\bar{R}}^{d}$ be the
product of its $d$ copies, so that $\mathbf{\bar{R}}^{d}$ itself is a
compactification of $\mathbf{R}^{d}$. (e.g. Dudley (2002), p.74.) We follow
the convention to set $\infty \cdot 0=0$ and $(-\infty )\cdot 0=0$. A
supremum and an infimum of a nonnegative map over an empty set are set to be
0 and $\infty $ respectively.

\section{Nondifferentiable Transforms of a Regular Parameter}

\subsection{Nondifferentiable Transforms}

First, we begin with conditions for $f$ and $g$ in (\ref{param2}).\bigskip

\noindent \textbf{Assumption 1:} (i) The map $g:\mathbf{R}^{d}\rightarrow 
\mathbf{R}$ is Lipschitz continuous, and satisfies the following.

(a) (Translation Equivariance) For each $c\in \mathbf{R}$ and $\mathbf{x}\in 
\mathbf{R}^{d}$, $g(\mathbf{x}+c)=g(\mathbf{x})+c.$

(b) (Scale Equivariance) For each $u\geq 0$ and $\mathbf{x}\in \mathbf{R}%
^{d},$ $g(u\mathbf{x})=ug(\mathbf{x}).$

(c) (Directional Derivatives) For each $\mathbf{z}\in \mathbf{R}^{d}$ and $%
\mathbf{x}\in \mathbf{R}^{d}$,%
\begin{equation*}
\tilde{g}(\mathbf{x};\mathbf{z})\equiv \lim_{t\downarrow 0}t^{-1}\left(
g\left( \mathbf{x}+t\mathbf{z}\right) -g\left( \mathbf{x}\right) \right)
\end{equation*}%
exists.

\noindent (ii) The map $f:\mathbf{\bar{R}}\rightarrow \mathbf{\bar{R}}$ is
Lipschitz continuous and non-constant on $\mathbf{R}$, and is continuously
differentiable except at a finite number of points in $\mathbf{R}$, with a
Lipschitz continuous derivative.\bigskip

Assumption 1(i) is the same as Assumption 1(i) of Song (2014) but the
requirement for $f$ is now substantially generalized by Assumption 1(ii). To
give a sense of the map $g$, consider the following examples.\bigskip

\noindent \textbf{Examples 1:} (a) $g(\mathbf{x})=\mathbf{s}^{\prime }%
\mathbf{x},$ where $\mathbf{s}\in S_{1}\equiv \{\mathbf{s}\in \mathbf{R}^{d}:%
\mathbf{s}^{\prime }\mathbf{1}_{d}=1\}$ and $\mathbf{1}_{d}$ is the $d$%
-dimensional vector of ones.

\noindent (b) $g(\mathbf{x})=\max (\mathbf{x})$ or $g(\mathbf{x})=\min (%
\mathbf{x})$.

\noindent (c) $g(\mathbf{x})=\max \{\min (\mathbf{x}_{1}),\mathbf{x}_{2}\}$,
where $\mathbf{x}_{1}$ and $\mathbf{x}_{2}$ are (possibly overlapping)
subvectors of $\mathbf{x}$.

\noindent (d) $g(\mathbf{x})=\max (\mathbf{x}_{1})+\max (\mathbf{x}_{2}),$ $%
g(\mathbf{x})=\min (\mathbf{x}_{1})+\min (\mathbf{x}_{2}),$ $g(\mathbf{x}%
)=\max (\mathbf{x}_{1})+\min (\mathbf{x}_{2}),$ or $g(\mathbf{x})=\max (%
\mathbf{x}_{1})+\mathbf{s}^{\prime }\mathbf{x}_{2}$ with $\mathbf{s}\in
S_{1}.$\bigskip

Thus the examples of the parameters $\theta $ in the form (\ref{param2}) are
as follows.\bigskip

\noindent \textbf{Examples 2:} (a) $\theta =|\max \{\beta _{1},\beta _{2}\}|$

\noindent (b) $\theta =|\beta _{1}-1|$.

\noindent (c) $\theta =\min \{\max \{\beta ,0\},1\}$.

\noindent (d) $\theta =\min \{|\beta |,1\}.$

\noindent (e) $\theta =\max \{\min \{\max \{\beta _{1},\beta _{2}\},1\},0\}.$%
\bigskip

The framework of Song (2014) requires that $f$ be a piecewise linear map
with a single kink point, and hence excludes the examples of (c)-(e). In
example (c), the parameter $\theta $ is $\beta $ censored at 0 and $1$. In
example (d), the parameter of interest is the absolute value of $|\beta |$
censored at 1.

One might ask whether the representation of parameter $\theta $ as a
composition map $f\circ g$ of $\mathbf{\beta }$ in (\ref{param2}) with $f$
and $g$ satisfying Assumption 1 is unique. Lemma 1 of Song (2014) gives an
affirmative answer. As we shall see later, the asymptotic risk bound
involves $g$ and the optimal estimators involve the maps $f$ and $g$. This
uniqueness result removes ambiguity that can potentially arise when $\theta $
has multiple equivalent representations with different maps $f$ and $g$.
When $d\geq 2$, the roles of $f$ and $g$ cannot be interchanged. When $d=1,$
Assumption 1(i) requires that $g(x)=x$. Hence in general the roles of $f$
and $g$ cannot be interchanged.

\subsection{Local Asymptotic Normality and Regularity of $\mathbf{\protect%
\beta }$}

We introduce briefly conditions for probabilities that identify $\mathbf{%
\beta }$, in a manner adapted from van der Vaart (1991) and van der Vaart
and Wellner (1996) (see Section 3.11, pp. 412-422.) Let $\mathcal{P}\equiv
\{P_{\alpha }:\alpha \in \mathcal{A}\}$ be a family of distributions on a
measurable space $(\mathcal{X},\mathcal{G})$ indexed by $\alpha \in \mathcal{%
A}$, where the set $\mathcal{A}$ is a subset of a Euclidean space or an
infinite dimensional space. We assume that we have i.i.d. draws $Y_{1},\cdot
\cdot \cdot ,Y_{n}$ from $P_{\alpha _{0}}\in \mathcal{P}$ so that $\mathbf{X}%
_{n}\equiv (Y_{1},\cdot \cdot \cdot ,Y_{n})$ is distributed as $P_{\alpha
_{0}}^{n}$. Let $\mathcal{P}(P_{\alpha _{0}})$ be the collection of maps $%
t\rightarrow P_{\alpha _{t}}\ $such that for some $h\in L_{2}(P_{\alpha
_{0}})$,%
\begin{equation}
\int \left\{ \frac{1}{t}\left( dP_{\alpha _{t}}^{1/2}-dP_{\alpha
_{0}}^{1/2}\right) -\frac{1}{2}hdP_{\alpha _{0}}^{1/2}\right\}
^{2}\rightarrow 0,\text{ as\ }n\rightarrow \infty .  \label{cv}
\end{equation}%
When this convergence holds, we say that $P_{\alpha _{t}}$ is \textit{%
differentiable in quadratic mean} to $P_{\alpha _{0}}$, call $h\in
L_{2}(P_{\alpha _{0}})$ a \textit{score function} associated with this
convergence, and call the set of all such $h$'s a \textit{tangent set},
denoting it by $T(P_{\alpha _{0}}).$ We assume that the tangent set is a
linear subspace of $L_{2}(P_{\alpha _{0}})$. Taking $\langle \cdot ,\cdot
\rangle $ to be the usual inner product in $L_{2}(P_{\alpha _{0}})$, we view 
$(H,\langle \cdot ,\cdot \rangle ),$ with $H\equiv T(P_{\alpha _{0}})$, as a
subspace of a separable Hilbert space. For each $h\in H,$ $n\in \mathbb{N}$,
and $\alpha _{h}\in \mathcal{A}$ such that $\alpha _{h}=0$ when $h=0$, let $%
P_{\alpha _{0}+\alpha _{h}/\sqrt{n}}$ be probabilities converging to $%
P_{\alpha _{0}}$ (as in (\ref{cv})) as $n\rightarrow \infty $ having $h$ as
its score. We simply write $P_{n,h}=P_{\alpha _{0}+\alpha _{h}/\sqrt{n}}^{n}$
and consider sequences of such probabilities $\{P_{n,h}\}_{n\geq 1}$ indexed
by $h\in H$. (See van der Vaart (1991) and van der Vaart and Wellner (1996),
Section 3.11 for details.) Differentiability in quadratic mean and i.i.d.
assumption imply local asymptotic normality (LAN): for any $h\in H$,%
\begin{equation*}
\log \frac{dP_{n,h}}{dP_{n,0}}=\zeta _{n}(h)-\frac{1}{2}\langle h,h\rangle ,
\end{equation*}%
where for each $h,h^{\prime }\in H$, 
\begin{equation*}
\lbrack \zeta _{n}(h),\zeta _{n}(h^{\prime })]\overset{d}{\rightarrow }%
[\zeta (h),\zeta (h^{\prime })],\text{\ under\ }\{P_{n,0}\},
\end{equation*}%
and $\zeta (\cdot )$ is a centered Gaussian process on $H$ with covariance
function $\mathbf{E}[\zeta (h_{1})\zeta (h_{2})]=\langle h_{1},h_{2}\rangle
. $ (See the proof of Lemma 3.10.11 of van der Vaart and Wellner (1996).)
Local asymptotic normality reduces the decision problem to one in which an
optimal decision is sought under a single Gaussian shift experiment $%
\mathcal{E}=(\mathcal{X},\mathcal{G},P_{h};h\in H),$ where $P_{h}$ is such
that $\log dP_{h}/dP_{0}=\zeta (h)-\frac{1}{2}\langle h,h\rangle .$

We assume that $\mathbf{\beta }$ is identified by $P_{\alpha }$, $\alpha \in 
\mathcal{A}$, and write $\mathbf{\beta }_{n}(h)=\mathbf{\beta }(P_{\alpha
_{0}+\alpha _{h}/\sqrt{n}}),$ regarding the parameter as an $\mathbf{R}^{d}$%
-valued map on $H$. As a sequence of maps on $H$, we assume that there
exists a continuous linear $\mathbf{R}^{d}$-valued map, $\mathbf{\dot{\beta}}
$, on $H$ such that for any $h\in H,$%
\begin{equation}
\sqrt{n}(\mathbf{\beta }_{n}(h)-\mathbf{\beta }_{n}(0))\rightarrow \mathbf{%
\dot{\beta}}(h)  \label{reg}
\end{equation}%
as $n\rightarrow \infty .$ In other words, $\mathbf{\beta }_{n}(h)\ $is 
\textit{regular} in the sense of van der Vaart \ and Wellner (1996, Section
3.11).

As is well known the functional $\mathbf{\dot{\beta}}$ determines the
efficiency bound for $\mathbf{\beta }$. Let $\mathbf{e}_{m}$ be a $d\times 1$
vector whose $m$-th entry is one and the other entries are zero. Certainly $%
\mathbf{e}_{m}^{\prime }\mathbf{\dot{\beta}}(\cdot )$ defines a continuous
linear functional on $H$, and hence there exists $\dot{\beta}_{m}^{\ast }\in 
\bar{H}$ such that $\mathbf{e}_{m}^{\prime }\mathbf{\dot{\beta}}(h)=\langle 
\dot{\beta}_{m}^{\ast },h\rangle ,$ $h\in H$. Then $||\dot{\beta}_{m}^{\ast
}||^{2}=\langle \dot{\beta}_{m}^{\ast },\dot{\beta}_{m}^{\ast }\rangle $
represents the asymptotic variance bound of the parameter $\beta _{m}=%
\mathbf{e}_{m}^{\prime }\mathbf{\beta }$. Let $\Sigma $ be a $d\times d$
matrix whose $(m,k)$-th entry is given by $\langle \dot{\beta}_{m}^{\ast },%
\dot{\beta}_{k}^{\ast }\rangle $. Throughout this paper, we assume that $%
\Sigma $ is invertible. The inverse of matrix $\Sigma $ is called the
semiparametric efficiency bound for $\mathbf{\beta }.$ (See Bickel,
Klaassen, Ritov and Wellner (1993) for ways to compute $\Sigma $.)

\section{Local Asymptotic Minimax Estimation}

\subsection{Loss Functions}

For a decision $d\in \mathbf{R}$ and the object of interest $\theta \in 
\mathbf{R}$, we consider the following form of a loss function:%
\begin{equation}
L\left( d,\theta \right) \equiv \tau (|d-\theta |),  \label{LB}
\end{equation}%
where $\tau :\mathbf{R}\rightarrow \mathbf{R}$ is a map that satisfies the
following assumption.\bigskip

\noindent \textbf{Assumption 2:} (i) $\tau (\cdot )$ is increasing on $%
[0,\infty )$, $\tau (0)=0$, and there exists $\bar{\tau}\in (0,\infty ]$ such%
$\ $that $\tau ^{-1}([0,y])$ is bounded in $[0,\infty )$ for all $0<y<\bar{%
\tau}$.

\noindent (ii) For each $M>0$, there exists $C_{M}>0$ such that for all $%
x,y\in \mathbf{R}$, 
\begin{equation}
|\tau _{M}(x)-\tau _{M}(y)|\ \leq \ C_{M}|x-y|,  \label{Lip}
\end{equation}%
where $\tau _{M}(\cdot )\equiv \min \{\tau (\cdot ),M\}$.\bigskip

The condition in (\ref{Lip}) allows unbounded loss functions. The class of
loss functions in this paper is mostly appropriate for the problem of
optimal estimation, but excludes some other types of decision problems. For
example, it excludes the hypothesis testing type loss function $\tau
(|d-\theta |)=1\{|d-\theta |>c\}$, $c\in \mathbf{R}$.

\subsection{Pointwise Local Asymptotic Minimax Theory}

First, we develop local asymptotic minimax theory for each \textit{fixed} $%
\alpha _{0}$, and call it \textit{pointwise local asymptotic minimax theory}%
, because the asymptotic approximation is pointwise at each $\alpha _{0}$.
Let%
\begin{equation*}
\theta _{n}(h)\equiv (f\circ g)(\mathbf{\beta }_{n}(h))
\end{equation*}%
and, given any estimator $\hat{\theta}$ which is a measurable function of $%
\mathbf{X}_{n}$, we define its \textit{local maximal risk}: for each $b\in
\lbrack 0,\infty ),$%
\begin{equation}
\mathcal{R}_{n,b}(\hat{\theta})\equiv \sup_{h\in H_{n,b}}\mathbf{E}_{h}\left[
\tau (|\sqrt{n}\{\hat{\theta}-\theta _{n}(h)\}|)\right] ,  \label{risk}
\end{equation}%
where $H_{n,b}\equiv \left\{ h\in H:\left\Vert \mathbf{\beta }_{n}(h)-%
\mathbf{\beta }_{n}(0)\right\Vert \leq b/\sqrt{n}\right\} $, and $\mathbf{E}%
_{h}$ denotes expectation under $P_{n,h}$.

Suppose that a Lipschitz continuous map $f:\mathbf{\bar{R}}\rightarrow 
\mathbf{\bar{R}}$ that satisfies Assumption 1(ii) is given. Let $\mathcal{Y}%
\subset \mathbf{R}$ be the set of differentiability points of $f$ in $%
\mathbf{R}$. By Assumption 1(ii), the set $\mathcal{Y}$ is dense in $\mathbf{%
R}$. We define for each $x\in \mathbf{R}$%
\begin{equation*}
\bar{f}^{\prime }(x)\equiv \ \underset{\varepsilon \downarrow 0}{\lim }%
\sup_{y\in \lbrack x-\varepsilon ,x+\varepsilon ]\cap \mathcal{Y}}\left\vert
f^{\prime }(y)\right\vert ,
\end{equation*}%
where $f^{\prime }(y)$ denotes the first order derivative of $f$ at $y$.
Certainly, the limit always exists by Assumption 1(ii), and hence $\bar{f}%
^{\prime }(x)$ is well defined for both differentiable and nondifferentiable
points. At a nondifferentiable point $x$, it is equivalent to define $\bar{f}%
^{\prime }(x)$ to be the maximum of the absolute left derivative and
absolute right derivative. One may consider various alternative concepts of
generalized derivatives (e.g. see Frank (1998)), but the definition $\bar{f}%
^{\prime }(x)$ is simple enough for our purpose. The following result is a
generalization of Theorem 1 in Song (2014).\bigskip

\noindent \textbf{Theorem 1:} \textit{Suppose that Assumptions 1-2\ hold.
Then for any sequence of estimators }$\hat{\theta}$,%
\begin{equation*}
\sup_{b\in \lbrack 0,\infty )}\underset{n\rightarrow \infty }{\text{liminf}}%
\ \mathcal{R}_{n,b}(\hat{\theta})\geq \inf_{c\in \mathbf{R}}B(c),
\end{equation*}%
\textit{where}%
\begin{equation*}
B(c)\equiv \sup_{\mathbf{r}\in \mathbf{R}^{d}}\mathbf{E}\left[ \tau \left( 
\bar{f}^{\prime }(g(\mathbf{\beta }_{0}))\left\vert \tilde{g}_{0}\left( Z+%
\mathbf{r}\right) -\tilde{g}_{0}(\mathbf{r})+c\right\vert \right) \right] .
\end{equation*}

The main feature of the local asymptotic risk bound in Theorem 1 is that it
involves infimum over a line, instead of infimum over an infinite
dimensional space. This convenient form is due to the same argument in Song
(2014) based on the purification result of Dvoretsky, Wald, and Wolfowitz
(1951) in zero sum games. This form is crucial for simulating the risk lower
bound when we construct a local asymptotic minimax estimator, as explained
below.

We consider an optimal estimator of $\theta $ that achieves the bound in
Theorem 1. The procedure here is adapted from the proposal by Song (2014).
Suppose that we are given a consistent estimator $\hat{\Sigma}$ of $\Sigma $
and a semiparametrically efficient estimator $\mathbf{\tilde{\beta}}$ of $%
\mathbf{\beta }$ which satisfy the following assumptions.\bigskip

\noindent \textbf{Assumption 3:} (i) For each $\varepsilon >0$, there exists 
$a>0$ such that 
\begin{equation*}
\text{limsup}_{n\rightarrow \infty }\text{sup}_{h\in H}\ P_{n,h}\{\sqrt{n}||%
\hat{\Sigma}-\Sigma ||>a\}<\varepsilon .
\end{equation*}

\noindent (ii) For each $t\in \mathbf{R}^{d}$, $\sup_{h\in H}\left\vert
P_{n,h}\{\sqrt{n}(\mathbf{\tilde{\beta}}-\mathbf{\beta }_{n}(h))\leq
t\}-P\{Z\leq t\}\right\vert \rightarrow 0$ as $n\rightarrow \infty $.\bigskip

Assumption 3 imposes $\sqrt{n}$-consistency of $\hat{\Sigma}$ and
convergence in distribution of $\sqrt{n}(\mathbf{\tilde{\beta}}-\mathbf{%
\beta }_{n}(h)),$ both uniform over $h\in H$. The uniform convergence can
often be verified through the central limit theorem uniform in $h\in H$.

For a fixed large $M_{1}>0,$ we define%
\begin{equation*}
\hat{\theta}_{mx}\equiv f\left( g(\mathbf{\tilde{\beta})}+\frac{\hat{c}%
_{M_{1}}}{\sqrt{n}}\right) ,
\end{equation*}%
where $\hat{c}_{M_{1}}$ is a bias adjustment term constructed from the
simulations of the risk lower bound in Theorem 1, as we explain now.

To simulate the risk lower bound in Theorem 1, we first draw $\{\mathbf{\xi }%
_{i}\}_{i=1}^{L}$ i.i.d. from $N(0,I_{d})$. Since $\tilde{g}_{0}(\cdot )$
depends on $\mathbf{\beta }_{0}$ that is unknown to the researcher, we first
construct a consistent estimator of $\tilde{g}_{0}(\cdot ).$ Take a sequence 
$\varepsilon _{n}\rightarrow 0$ such that $\sqrt{n}\varepsilon
_{n}\rightarrow \infty $ as $n\rightarrow \infty $. Examples of $\varepsilon
_{n}$ are $\varepsilon _{n}=n^{-1/3}$ or $\varepsilon _{n}=n^{-1/2}\log n$.
Let%
\begin{equation*}
\hat{g}_{n}(\mathbf{z})\equiv g\left( \mathbf{z+}\varepsilon _{n}^{-1}(%
\mathbf{\tilde{\beta}}-g(\mathbf{\tilde{\beta}}))\right) .
\end{equation*}%
Then it is not hard to see that $\hat{g}_{n}(\mathbf{z})$ is consistent for $%
\tilde{g}_{0}(\mathbf{z})$. Define%
\begin{equation*}
\hat{a}_{n}\equiv \sup_{x\in \lbrack g(\mathbf{\tilde{\beta})}-\varepsilon
_{n},g(\mathbf{\tilde{\beta})}+\varepsilon _{n}]\cap \mathcal{Y}}\left\vert
f^{\prime }\left( x\right) \right\vert .
\end{equation*}%
Let%
\begin{equation}
\hat{B}_{M_{1}}(c)\equiv \sup_{\mathbf{r}\in \lbrack -M_{1},M_{1}]^{d}}\frac{%
1}{L}\sum_{i=1}^{L}\tau _{M_{1}}\left( \hat{a}_{n}\left\vert \hat{g}_{n}(%
\hat{\Sigma}^{1/2}\mathbf{\xi }_{i}+\mathbf{r})-\hat{g}_{n}(\mathbf{r}%
)+c\right\vert \right) \text{.}  \label{BM}
\end{equation}%
Then we define%
\begin{equation}
\hat{c}_{M_{1}}\equiv \frac{1}{2}\left\{ \sup \hat{E}_{M_{1}}+\inf \hat{E}%
_{M_{1}}\right\} ,  \label{cmx}
\end{equation}%
where, with $\eta _{n,L}\rightarrow 0$, $\eta _{n,L}(\sqrt{L}+\varepsilon
_{n}\sqrt{n}+\varepsilon _{n}^{-1})\rightarrow \infty $ as$\ n,L\rightarrow
\infty $,%
\begin{equation*}
\hat{E}_{M_{1}}\equiv \left\{ c\in \lbrack -M_{1},M_{1}]:\hat{B}%
_{M_{1}}(c)\leq \inf_{c_{1}\in \lbrack -M_{1},M_{1}]}\hat{B}%
_{M_{1}}(c_{1})+\eta _{n,L}\right\} .
\end{equation*}

The following theorem affirms that $\hat{\theta}_{mx}$ is local asymptotic
minimax for $\theta =g(\mathbf{\beta })$. (For technical facility, we follow
a suggestion by Strasser (1985) (p.440) and consider a truncated loss: $\tau
_{M}(\cdot )=\min \{\tau (\cdot ),M\}$ for large $M.$)\bigskip

\noindent \textbf{Theorem 2:} \textit{Suppose that Assumptions 1-3\ hold.
Then, for any }$M>0$\textit{\ and any }$M_{1}\geq M$\textit{\ that
constitutes constant }$\hat{c}_{M_{1}},$%
\begin{equation*}
\sup_{b\in \lbrack 0,\infty )}\underset{n\rightarrow \infty }{\text{limsup}}%
\ \mathcal{R}_{n,b,M}(\hat{\theta}_{mx})\leq \inf_{c\in \mathbf{R}}B(c),
\end{equation*}%
where $\mathcal{R}_{n,b,M}(\hat{\theta}_{mx})$\textit{\ coincides with }$%
\mathcal{R}_{n,b}(\hat{\theta}_{mx})$ \textit{with }$\tau (\cdot )$\textit{\
replaced by }$\min \{\tau (\cdot ),M\}$.\bigskip

Therefore, the risk lower bound 
\begin{equation*}
\inf_{c\in \mathbf{R}}B(c)
\end{equation*}%
in Theorem 1 is sharp. We call it the \textit{local asymptotic minimax risk}
in this paper.

When $\tau (x)=|x|^{p}$, for some $p\geq 1$, the minimizer of $B(c)$ does
not depend on the shape of $f$. Hence, in constructing $\hat{B}_{M_{1}}(c)$
in (\ref{BM}), it suffices to take $\hat{a}_{n}=1$.

When $\theta =g(\mathbf{\beta })$ is a regular parameter, taking the form of 
$g(\mathbf{\beta })=\mathbf{s}^{\prime }\mathbf{\beta }$ with $\mathbf{s}\in
S_{1}$, the local asymptotic minimax risk bound becomes%
\begin{equation*}
\inf_{c\in \mathbf{R}}\mathbf{E}\left[ \tau \left( \bar{f}^{\prime }(\mathbf{%
s}^{\prime }\mathbf{\beta }_{0}))|\mathbf{s}^{\prime }Z+c|\right) \right] =%
\mathbf{E}\left[ \tau \left( \bar{f}^{\prime }(\mathbf{s}^{\prime }\mathbf{%
\beta }_{0})|\mathbf{s}^{\prime }Z|\right) \right] ,
\end{equation*}%
where the equality above follows by Anderson's Lemma. In this case, it
suffices to set $\hat{c}_{M_{1}}^{\ast }=0$, for the infimum over $c\in 
\mathbf{R}$ is achieved at $c=0$. This is true regardless of whether $f$ is
symmetric around zero or not. Hence the minimax decision becomes simply%
\begin{equation}
\hat{\theta}_{\text{$mx$}}=f(\mathbf{\tilde{\beta}}^{\prime }\mathbf{s}).
\label{sol}
\end{equation}%
This has the following consequences.\bigskip

\noindent \textbf{Example 5:} (a) When $\theta =\mathbf{\beta }^{\prime }%
\mathbf{s}$ for a known vector $\mathbf{s}\in S_{1}$, $\hat{\theta}_{\text{$%
mx$}}=\mathbf{\tilde{\beta}}^{\prime }\mathbf{s}$. Therefore, the decision
in (\ref{sol}) reduces to the well-known semiparametric efficient estimator
of $\mathbf{\beta }^{\prime }\mathbf{s}$.

\noindent (b) When $\theta =\max \{a\cdot \mathbf{\beta }^{\prime }\mathbf{s}%
+b,0\}$ for a known vector $\mathbf{s}\in S_{1}$ and known constants $a,b\in 
\mathbf{R}$, $\hat{\theta}_{\text{$mx$}}=\max \{a\cdot \mathbf{\tilde{\beta}}%
^{\prime }\mathbf{s}+b,0\}.$

\noindent (c) When $\theta =|\beta |$ for a scalar parameter $\beta $, $\hat{%
\theta}_{\text{$mx$}}=|\hat{\beta}|.$ This decision is analogous to
Blumenthal and Cohen (1968).

\noindent (d) When $\theta =\max \{\beta _{1}+\beta _{2}-1,0\}$, $\hat{\theta%
}_{mx}=\max \{\hat{\beta}_{1}+\hat{\beta}_{2}-1,0\}$. $\blacksquare $\bigskip

The examples of (b)-(d)\ involve nondifferentiable transform $f$, and hence $%
\hat{\theta}_{\text{$mx$}}=f(\mathbf{\tilde{\beta}}^{\prime }\mathbf{s)}$ as
an estimator of $\theta $ is asymptotically biased in these examples.
Nevertheless, the plug-in estimator $\hat{\theta}_{\text{$mx$}}$ that does
not involve any bias-reduction is local asymptotic minimax.

\section{Discontinuity in the Local Asymptotic Minimax Risk}

\subsection{Local Robustification}

The local asymptotic minimax risk $\inf_{c\in \mathbf{R}}B(c)$ depends on $%
\mathbf{\beta }_{0}$ discontinuously in general. This is easily seen from
the form of $B(c)$, for we may have $\bar{f}^{\prime }(x)$ and $\bar{f}%
^{\prime }(y)$ stay apart, even as $x$ and $y$ get closer to each other.%
\footnote{%
While nondifferentiability of $f$ yields discontinuity in the minimax risk
in most cases, there are counterexamples. For example, when $f(a)=|a|$, we
have $\bar{f}^{\prime }(a)=1$ for all $a\in \mathbf{R,}$ and if further $g(%
\mathbf{\beta })=\mathbf{s}^{\prime }\mathbf{\beta }$, the minimax risk is
continuous in $\alpha _{0}$, although $f$ is nondifferentiable at $0$.} This
discontinuity may imply that the local asymptotic minimax approach may serve
as a poor approximation of a finite sample risk bound.

We consider an alternative approach of optimality that is robustified
against a local perturbation of $\mathbf{\beta }_{0}=\mathbf{\beta }%
(P_{\alpha _{0}})$ (i.e. of $\alpha _{0}\in \mathcal{A}$). Note that
Ibragimov and Khas'minski (1986) pointed out the desirability of local
robustification with respect to such an initial parameter $\mathbf{\beta }%
_{0}$.

For each $\mathbf{y\in R}^{d}$ and a positive sequence $\varepsilon
_{n}\downarrow 0$, define%
\begin{equation*}
\mathcal{A}(\alpha _{0};\varepsilon _{n})\equiv \left\{ \alpha \in \mathcal{A%
}:||\mathbf{\beta }(P_{\alpha })-\mathbf{\beta }(P_{\alpha _{0}})||\leq
\varepsilon _{n}\right\} .
\end{equation*}%
The set $\mathcal{A}(\alpha _{0};\varepsilon _{n})$ is the collection of $%
\alpha $'s such that the regular parameter vectors $\mathbf{\beta }%
(P_{\alpha })$ and $\mathbf{\beta }(P_{\alpha _{0}})$ are close to each
other. Then define the local maximal risk under local robustification: for
each $b\in \lbrack 0,\infty ),$ and a positive sequence $\varepsilon
_{n}\downarrow 0,$%
\begin{equation*}
\mathcal{R}_{n,b}(\hat{\theta};\varepsilon _{n})\equiv \sup_{\alpha _{1}\in 
\mathcal{A}(\alpha _{0};\varepsilon _{n})}\sup_{h\in H_{n,b}(\alpha _{1})}%
\mathbf{E}_{h,\alpha _{1}}\left[ \tau (|\sqrt{n}\{\hat{\theta}-\theta
_{n}(h)\}|)\right] ,
\end{equation*}%
where $\mathbf{E}_{h,\alpha _{1}}$ denotes the expectation under $%
P_{n,h,\alpha _{1}}\equiv P_{\alpha _{1}+\alpha _{h}/\sqrt{n}}$ and 
\begin{equation*}
H_{n,b}(\alpha _{1})\equiv \left\{ h\in H:\left\Vert \mathbf{\beta }_{n}(h)-%
\mathbf{\beta }(P_{\alpha _{1}})\right\Vert \leq b/\sqrt{n}\right\} .
\end{equation*}%
Then certainly by Theorem 1, we have for any sequence of estimators $\hat{%
\theta}$,%
\begin{equation*}
\sup_{b\in \lbrack 0,\infty )}\underset{n\rightarrow \infty }{\text{liminf}}%
\ \mathcal{R}_{n,b}(\hat{\theta};\varepsilon _{n})\geq \inf_{c\in \mathbf{R}%
}B(c).
\end{equation*}%
The main question is whether this lower bound is sharp. For this, we show
that the optimal estimator $\hat{\theta}_{mx}$ continues to achieve this
lower bound, when Assumption 3 is strengthened as follows.\bigskip

\noindent \textbf{Assumption 3':} (i) There exists $M>0$ such that 
\begin{equation*}
\text{limsup}_{n\rightarrow \infty }\text{sup}_{\alpha \in \mathcal{A}}\text{%
sup}_{h\in H}P_{n,h,\alpha }\{\sqrt{n}||\hat{\Sigma}-\Sigma
||>M\}<\varepsilon .
\end{equation*}

\noindent (ii) For each $t\in \mathbf{R}^{d}$, 
\begin{equation*}
\sup_{\alpha \in \mathcal{A}}\sup_{h\in H}\left\vert P_{n,h,\alpha }\{\sqrt{n%
}(\mathbf{\tilde{\beta}}-\mathbf{\beta }_{n}(h))\leq t\}-P\{Z\leq
t\}\right\vert \rightarrow 0,
\end{equation*}%
as $n\rightarrow \infty $.\bigskip

Assumption 3' strengthens the uniformity in convergence in Assumption 3 to
that over $\alpha \in \mathcal{A}$. In many cases, it is not hard to verify
this condition.\bigskip

\noindent \textbf{Theorem 3:} \textit{Suppose that Assumptions 1-2 and 3'\
hold. Then for each }$\varepsilon _{n}\downarrow 0$\textit{\ such that }$%
\varepsilon _{n}\sqrt{n}\rightarrow \infty ,$\textit{\ as }$n\rightarrow
\infty $\textit{, for any sequence of estimators }$\hat{\theta}$\textit{,
and for any }$M_{1}>M$ \textit{such that} $M_{1}$ \textit{constitutes} $\hat{%
c}_{M_{1}}$,%
\begin{equation*}
\lim_{M\uparrow \infty }\sup_{b\in \lbrack 0,\infty )}\underset{n\rightarrow
\infty }{\text{limsup}}\ \mathcal{R}_{n,b,M}(\hat{\theta}_{mx};\varepsilon
_{n})\leq \inf_{c\in \mathbf{R}}B(c).
\end{equation*}

The result of Theorem 3 shows that the local asymptotic minimax risk 
\begin{equation}
\inf_{c\in \mathbf{R}}B(c)  \label{mr}
\end{equation}%
remains unchanged, even when we locally robustify the maximal risk, and that
the estimator $\hat{\theta}_{mx}$ continues to satisfy the local asymptotic
minimaxity after local robustification. At the same time, local
robustification does not resolve the issue of discontinuity in the local
asymptotic minimax risk.

Given Theorems 1-3, we find that fixing $\bar{\delta}>0$, and considering $%
\mathcal{A}(\alpha _{0};\bar{\delta}/\sqrt{n})$ instead of $\mathcal{A}%
(\alpha _{0};\varepsilon _{n})$ will not change the result. The local
asymptotic minimax risk does not depend on this choice of $\bar{\delta}$.
This is because the local asymptotic minimax risk remains the same either we
take $\mathcal{R}_{n,b}(\hat{\theta};0)$ (Theorems 1-2) or we take $\mathcal{%
R}_{n,b}(\hat{\theta};\varepsilon _{n})$ as our local maximal risk.

\subsection{Discussion}

To understand the result of Theorem 3, let us consider localization with a
fixed Pitman direction. First, let $\alpha _{0}\in \mathcal{A}$ be as before
such that $\mathbf{\beta }(P_{\alpha _{0}})=\mathbf{\beta }_{0}$, and
consider 
\begin{equation}
\alpha _{n}(\delta )=\alpha _{0}+\frac{\delta }{\sqrt{n}},  \label{lp}
\end{equation}%
where $\delta \in \mathcal{A}$ is a Pitman direction. We assume that $%
\{P_{\alpha _{n}(\delta )}\}_{n=1}^{\infty }$ is quadratic mean
differentiable at $P_{\alpha _{0}}$ with a score $h_{\delta }\in H$. Then,
we have 
\begin{eqnarray*}
\mathbf{\beta }(P_{\alpha _{n}(\delta )}) &=&\mathbf{\beta }(P_{\alpha
_{0}})+\mathbf{\beta }(P_{\alpha _{0}+\delta /\sqrt{n}})-\mathbf{\beta }%
(P_{\alpha _{0}}) \\
&=&\mathbf{\beta }(P_{\alpha _{0}})+\frac{\mathbf{\dot{\beta}}(h_{\delta })}{%
\sqrt{n}}+o(n^{-1/2}),
\end{eqnarray*}%
by the regularity of $\mathbf{\beta }$. In other words, the Pitman direction 
$\delta $ for $\alpha _{n}(\delta )$ is now translated into the Pitman
direction $\mathbf{\dot{\beta}}(h_{\delta })$ for $\mathbf{\beta }(P_{\alpha
_{n}(\delta )})$.

Recall that the local asymptotic minimax risk arises as a consequence of
robustification against all the scores $h\in H$ at $P_{\alpha _{0}}$. Since $%
\Sigma $ is invertible, the range of $\mathbf{\dot{\beta}}$ (when $\mathbf{%
\dot{\beta}}$ is extended to a completion $\bar{H}$ of $H$) is equal to $%
\mathbf{R}^{d}$, i.e, for any $\mathbf{r}\in \mathbf{R}^{d},$ there exists $%
h\in \bar{H}$ such that $\mathbf{\dot{\beta}}(h)=\mathbf{r}$. Therefore,
robustification against all the Pitman directions such that $\{P_{\alpha
_{n}(\delta )}\}_{n=1}^{\infty }$ is quadratic mean differentiable at $%
P_{\alpha _{0}}$ is equivalent to robustification against all the $\sqrt{n}$%
-converging Pitman deviations from $\mathbf{\beta }_{0}$. Thus the local
robustification against Pitman deviations from $\mathbf{\beta }$ is already
incorporated in the results of local asymptotic minimax risk in Theorems 1
and 2.\footnote{%
In fact, Theorem 3 is stronger than this, because for each $\bar{\delta}\in
(0,\infty ),$ 
\begin{equation*}
\mathcal{A}(\alpha _{0};\bar{\delta}/\sqrt{n})\subset \mathcal{A}(\alpha
_{0};\varepsilon _{n}),
\end{equation*}%
from some large $n$ on.} This is why local robustification around $\mathbf{%
\beta }_{0}$ does not alter the results.

One might suggest considering a single Pitman direction $\delta $ and
focusing on a sequence of probabilities $\{P_{\alpha _{n}(\delta
)}\}_{n=1}^{\infty }$, derive the local asymptotic minimax risk in a way
that depends on $\delta $, and see if the risk continuously depends on $%
\delta $. This approach is analogous to many other approaches used to deal
with discontinuity of asymptotic distributions such as Pitman local
asymptotic power analysis, local-to-unity models, and weak identification.
However, such an approach in this context counters to the basic motivation
of the local asymptotic minimax approach, because restricting attention to a
single sequence of probabilities fails to robustify the decision problem
properly against local perturbations of the underlying probability and hence
fails to exclude superefficient estimators.\footnote{%
A still alternative way is an approach of global robustification, where one
robustifies against all $\alpha _{0}$'s in $\mathcal{A}$. The problem with
this approach is that the minimax decision problem often becomes trivial,
with the minimax risk being infinity. Such a trivial case arises, for
example, when $\sup_{x\in \mathbf{R}}\bar{f}^{\prime }(x)=\infty $. This is
the case when $f(x)=x^{2}$ for example.}

\section{Conclusion}

This paper focuses on the problem of optimal estimation for a parameter that
is a nondifferentiable transform of a regular parameter. First, this paper
extends the results of Song (2014) allowing for a more general class of
nondifferentiable transforms. Second, this paper investigates the issue of
discontinuity in local asymptotic minimax risk, and considers the approach
of local robustification of the true probability. As it turns out, the local
robustification does not alter the local asymptotic minimax risk. This means
that the optimal estimator remains optimal under this additional dimension
of local robustification. On the other hand, it also means that the
discontinuity in the minimax risk is not resolved by the local
robustification. Hence, there still remains the question of whether local
asymptotic minimax theory gives a good approximation of a finite sample
decision problem when the parameter is nondifferentiable. A\ full
investigation of this issue is relegated to a future research.

\section{Appendix: Mathematical Proofs}

\noindent \textbf{Proof of Theorem 1:} As in the proof of Lemma 3 of Song
(2014), we begin by fixing $\mathbf{r}\in \mathbf{R}^{d}$ so that for some $%
h^{\prime }\in \overline{H}$, $\mathbf{r}=\mathbf{\dot{\beta}}(h^{\prime })$%
. (Existence of such $h^{\prime }\in \overline{H}$ for each $\mathbf{r\in R}%
^{d}$ follows from the condition that $\Sigma $ is invertible.) Also fix $%
h\in H$ such that $\langle h,h^{\prime }\rangle =0$. We write%
\begin{eqnarray*}
&&\sqrt{n}\{\hat{\theta}-f(g(\mathbf{\beta }_{n}(h+h^{\prime })))\} \\
&=&\sqrt{n}\{\hat{\theta}-f(g(\mathbf{\beta }_{n}(h^{\prime })))\}-\sqrt{n}%
\left\{ f(g(\mathbf{\beta }_{n}(h+h^{\prime })))-f(g(\mathbf{\beta }%
_{n}(h^{\prime })))\right\} .
\end{eqnarray*}

Suppose first that $f$ is continuously differentiable at $g(\mathbf{\beta }%
_{0})$. Then by Assumption 1(ii), we have $x_{1}<x_{2}$ such that $x_{1}<g(%
\mathbf{\beta }_{0})<x_{2}$, where $f$ is continuously differentiable on $%
[x_{1},x_{2}]$. Furthermore, by regularity of $\mathbf{\beta }$ and
Lipschitz continuity of $g$, we have 
\begin{equation}
g(\mathbf{\beta }_{n}(h+h^{\prime }))=g(\mathbf{\beta }_{n}(h+h^{\prime })-%
\mathbf{\beta }_{0}+\mathbf{\beta }_{0})\rightarrow g(\mathbf{\beta }_{0}),
\label{cv7}
\end{equation}%
as $n\rightarrow \infty $. Therefore, we note that from some large $n$ on,
by the mean value theorem, 
\begin{eqnarray*}
&&\sqrt{n}\left\{ f(g(\mathbf{\beta }_{n}(h+h^{\prime })))-f(g(\mathbf{\beta 
}_{0}))\right\} \\
&=&\sqrt{n}f^{\prime }(a_{n}(h,h^{\prime }))\left\{ g(\mathbf{\beta }%
_{n}(h+h^{\prime }))-g(\mathbf{\beta }_{0})\right\} ,
\end{eqnarray*}%
where $a_{n}(h,h^{\prime })\equiv t_{n}\{g(\mathbf{\beta }_{n}(h+h^{\prime
})-g(\mathbf{\beta }_{0})\}+g(\mathbf{\beta }_{0})$ for some $t_{n}\in
\lbrack 0,1]$. From (\ref{cv7}), we have 
\begin{equation*}
a_{n}(h,h^{\prime })\rightarrow g(\mathbf{\beta }_{0})\text{, as }%
n\rightarrow \infty \text{.}
\end{equation*}%
From (A.10) of Song (2014) on page 149, we also find that%
\begin{equation}
\sqrt{n}\{g(\mathbf{\beta }_{n}(h+h^{\prime }))-g(\mathbf{\beta }%
_{n}(h^{\prime }))\}=\tilde{g}_{0}(\mathbf{\dot{\beta}}(h)+\mathbf{r})-%
\tilde{g}_{0}(\mathbf{r})+o(1).  \label{cv6}
\end{equation}%
Since $f^{\prime }$ is continuous at $g(\mathbf{\beta }_{0})\ $(Assumption
1(ii)) and $\mathbf{\dot{\beta}}$ is bounded, we combine these results to
deduce that%
\begin{eqnarray}
&&\sqrt{n}\left\{ f(g(\mathbf{\beta }_{n}(h+h^{\prime })))-f(g(\mathbf{\beta 
}_{n}(h^{\prime })))\right\}  \label{cv8} \\
&\rightarrow &f^{\prime }(g(\mathbf{\beta }_{0}))\left( \tilde{g}_{0}(%
\mathbf{\dot{\beta}}(h)+\mathbf{r})-\tilde{g}_{0}(\mathbf{r})\right) , 
\notag
\end{eqnarray}%
as $n\rightarrow \infty $.

For any sequence of estimators $\hat{\theta}$, the sequence $\{\hat{\theta}%
\}_{n\geq 1}$ is uniformly tight in $\mathbf{\bar{R}}$, and hence by using
the LAN property and (\ref{cv8}), applying Prohorov's Theorem, we find that
for each subsequence of $\{n\}$, there exists a further subsequence $%
\{n^{\prime }\}$ such that under $\{P_{n^{\prime },h^{\prime }}\},$ 
\begin{equation*}
\left[ 
\begin{array}{c}
\sqrt{n^{\prime }}\left\{ \hat{\theta}-g(\mathbf{\beta }_{n^{\prime
}}(h+h^{\prime }))\right\} \\ 
\log dP_{n^{\prime },h+h^{\prime }}/dP_{n^{\prime },h^{\prime }}%
\end{array}%
\right] \overset{d}{\rightarrow }\left[ 
\begin{array}{c}
V-f^{\prime }(g(\mathbf{\beta }_{0}))\left( \tilde{g}_{0}(\mathbf{\dot{\beta}%
}(h)+\mathbf{r})-\tilde{g}_{0}(\mathbf{r})\right) \\ 
\zeta (h)-\frac{1}{2}\langle h,h\rangle%
\end{array}%
\right] ,
\end{equation*}%
where $V\in \mathbf{\bar{R}}$ is a random variable having a potentially
deficient distribution. The rest of the proof can be proceeded precisely as
in the proofs of Lemma 3 and Theorem 1 of Song (2014).

Second, suppose that $f$ is not continuously differentiable at $g(\mathbf{%
\beta }_{0})$. Since there is a finite number of nondifferentiability points
for $f$, we have $x_{1},x_{2}\in \mathbf{R}$ such that $x_{1}<g(\mathbf{%
\beta }_{0})<x_{2}$, and $f$ is continuously differentiable on $[x_{1},g(%
\mathbf{\beta }_{0}))$ and $(g(\mathbf{\beta }_{0}),x_{2}]$ by Assumption
1(ii).

As previously, we choose arbitrary $\mathbf{r}\in \mathbf{R}^{d}$ so that
for some $h^{\prime }\in H$, $\mathbf{r=\dot{\beta}(}h^{\prime }\mathbf{)}$,
and fix $b/2\geq ||h^{\prime }||\cdot ||\mathbf{\dot{\beta}}^{\ast }||$.
Now, define%
\begin{eqnarray*}
H_{n,b,1}^{\ast } &\equiv &\{h\in H_{n,b}^{\ast }:R_{n}(h)\geq 0\},\text{ and%
} \\
H_{n,b,2}^{\ast } &\equiv &\{h\in H_{n,b}^{\ast }:R_{n}(h)\leq 0\},
\end{eqnarray*}%
where $H_{n,b}^{\ast }\equiv \{h\in H_{n,b}:\langle h,h^{\prime }\rangle
=0\} $ and $R_{n}(h)\equiv g(\mathbf{\beta }_{n}(h+h^{\prime }))-g(\mathbf{%
\beta }_{0})$, and observe that for all $h\in H,$%
\begin{equation}
\left\vert R_{n}(h)\right\vert \rightarrow 0,  \label{cv2}
\end{equation}%
as $n\rightarrow \infty $. (This is (A.31) of Song (2014). See the arguments
for details.) Thus we have for each $h\in H$, as $n\rightarrow \infty $,%
\begin{eqnarray}
1\left\{ h\in H_{n,b,1}^{\ast }\right\} &\rightarrow &1\left\{ h\in
H_{b}\right\} \text{ and}  \label{cv3} \\
1\left\{ h\in H_{n,b,2}^{\ast }\right\} &\rightarrow &1\left\{ h\in
H_{b}\right\} ,  \notag
\end{eqnarray}%
where $H_{b}\equiv \{h\in H:||\mathbf{\dot{\beta}}(h)||\leq b\}$.

Note that 
\begin{eqnarray}
&&\underset{b\rightarrow \infty }{\text{lim}}\ \underset{n\rightarrow \infty 
}{\text{liminf}}\sup_{h\in H_{n,b}}\mathbf{E}_{h}\left[ \tau _{M}\left( |%
\sqrt{n}\{\hat{\theta}-f(g(\mathbf{\beta }_{n}(h)))|\right) \right]
\label{bd7} \\
&\geq &\max_{l=1,2}\underset{b\rightarrow \infty }{\text{lim}}\ \underset{%
n\rightarrow \infty }{\text{liminf}}\sup_{h\in H_{n,b/2,l}^{\ast }}\mathbf{E}%
_{h}\left[ \tau _{M}\left( |\sqrt{n}\{\hat{\theta}-f\left( g(\mathbf{\beta }%
_{n}(h))\right) |\right) \right] .  \notag
\end{eqnarray}%
Due to the liminf and supremum over $h\in H_{n,b/2,l}^{\ast }$ (where the
supremum over an empty set of a nonnegative function is taken to be zero),
it suffices to focus on $h\in H$ such that $h\in H_{n,b/2,1}^{\ast }$ or $%
h\in H_{n,b/2,2}^{\ast }$, eventually from some large $n$ on.

For each $h\in H_{n,b,1}^{\ast }$, we have 
\begin{equation*}
g(\mathbf{\beta }_{0})\leq g(\mathbf{\beta }_{n}(h+h^{\prime })),
\end{equation*}%
and by the mean-value theorem, 
\begin{eqnarray*}
&&\sqrt{n}\left\{ f(g(\mathbf{\beta }_{n}(h+h^{\prime })))-f(g(\mathbf{\beta 
}_{0}))\right\} \\
&=&f_{+}^{\prime }(g(\mathbf{\beta }_{0})+t_{n})\sqrt{n}\left( g(\mathbf{%
\beta }_{n}(h+h^{\prime }))-g(\mathbf{\beta }_{0})\right) ,
\end{eqnarray*}%
where $t_{n}\geq 0$ and $t_{n}\leq g(\mathbf{\beta }_{n}(h+h^{\prime }))-g(%
\mathbf{\beta }_{0})$ and $f_{+}^{\prime }(x)$ denotes the right derivative
of $f$ at $x$. By using (\ref{cv6}) and the Lipschitz continuity of $%
f^{\prime }$ on $(g(\mathbf{\beta }_{0}),x_{2}]$, for any $h\in H,$ such
that $h\in H_{n,b/2,1}^{\ast }$ eventually, we have%
\begin{equation*}
\sqrt{n}\left\{ f(g(\mathbf{\beta }_{n}(h+h^{\prime })))-f(g(\mathbf{\beta }%
_{n}(h^{\prime })))\right\} \rightarrow f_{+}^{\prime }(g(\mathbf{\beta }%
_{0}))\left( \tilde{g}_{0}(\mathbf{\dot{\beta}}(h)+\mathbf{r})-\tilde{g}_{0}(%
\mathbf{r})\right) ,
\end{equation*}%
as $n\rightarrow \infty $.

Therefore, for any $h\in H,$ such that $h\in H_{n,b/2,1}^{\ast }$
eventually, and for each subsequence of $\{n\}$, there exists a further
subsequence $\{n^{\prime }\}$ such that under $\{P_{n^{\prime },h^{\prime
}}\},$%
\begin{equation*}
\left[ 
\begin{array}{c}
\sqrt{n^{\prime }}\left\{ \hat{\theta}-f\left( g(\mathbf{\beta }_{n^{\prime
}}(h+h^{\prime }))\right) \right\} \\ 
\log dP_{n^{\prime },h+h^{\prime }}/dP_{n^{\prime },h^{\prime }}%
\end{array}%
\right] \overset{d}{\rightarrow }\left[ 
\begin{array}{c}
V_{+}-f_{+}^{\prime }(g(\mathbf{\beta }_{0}))\left( \tilde{g}_{0}(\mathbf{%
\dot{\beta}}(h)+\mathbf{r})-\tilde{g}_{0}(\mathbf{r})\right) \\ 
\zeta (h)-\frac{1}{2}\langle h,h\rangle%
\end{array}%
\right] ,
\end{equation*}%
where $V_{+}\in \mathbf{\bar{R}}$ is a random variable having a potentially
deficient distribution. Using (\ref{cv3}) and following the same arguments
as in the proofs of Lemma 3 and Theorem 1 of Song (2014), we deduce that 
\begin{eqnarray}
&&\underset{b\rightarrow \infty }{\text{lim}}\ \underset{n\rightarrow \infty 
}{\text{liminf}}\sup_{h\in H_{n,b/2,1}^{\ast }}\mathbf{E}_{h}\left[ \tau
_{M}\left( |\sqrt{n}\{\hat{\theta}-f\left( g(\mathbf{\beta }_{n}(h))\right)
|\right) \right]  \label{st1} \\
&\geq &\inf_{c\in \mathbf{R}}\sup_{\mathbf{r}\in \mathbf{R}^{d}}\mathbf{E}%
\left[ \tau \left( \left\vert f_{+}^{\prime }(g(\mathbf{\beta }%
_{0}))\right\vert \left\vert \tilde{g}_{0}\left( Z+\mathbf{r}\right) -\tilde{%
g}_{0}(\mathbf{r})+c\right\vert \right) \right] .  \notag
\end{eqnarray}

Similarly, for any $h\in H,$ such that $h\in H_{n,b/2,2}^{\ast }$
eventually, we have%
\begin{equation*}
\sqrt{n}\left\{ f(g(\mathbf{\beta }_{n}(h_{n}+h^{\prime })))-f(g(\mathbf{%
\beta }_{n}(h^{\prime })))\right\} \rightarrow f_{-}^{\prime }(g(\mathbf{%
\beta }_{0}))\left( \tilde{g}_{0}(\mathbf{\dot{\beta}}(h)+\mathbf{r})-\tilde{%
g}_{0}(\mathbf{r})\right) ,
\end{equation*}%
as $n\rightarrow \infty $, where $f_{-}^{\prime }(x)$ denotes the left
derivative of $f$ at $x$. Hence similarly as before, for any $h\in H,$ such
that $h\in H_{n,b/2,2}^{\ast }$ eventually, and for each subsequence of $%
\{n\}$, there exists a further subsequence $\{n^{\prime }\}$ such that under 
$\{P_{n^{\prime },h^{\prime }}\},$%
\begin{equation*}
\left[ 
\begin{array}{c}
\sqrt{n^{\prime }}\left\{ \hat{\theta}-f\left( g(\mathbf{\beta }_{n^{\prime
}}(h+h^{\prime }))\right) \right\} \\ 
\log dP_{n^{\prime },h+h^{\prime }}/dP_{n^{\prime },h^{\prime }}%
\end{array}%
\right] \overset{d}{\rightarrow }\left[ 
\begin{array}{c}
V_{-}-f_{-}^{\prime }(g(\mathbf{\beta }_{0}))\left( \tilde{g}_{0}(\mathbf{%
\dot{\beta}}(h)+\mathbf{r})-\tilde{g}_{0}(\mathbf{r})\right) \\ 
\zeta (h)-\frac{1}{2}\langle h,h\rangle%
\end{array}%
\right] ,
\end{equation*}%
where $V_{-}\in \mathbf{\bar{R}}$ is a random variable having a potentially
deficient distribution. Using (\ref{cv3}) and following the same arguments
as in the proofs of Lemma 3 and Theorem 1 of Song (2014), we deduce that 
\begin{eqnarray}
&&\underset{b\rightarrow \infty }{\text{lim}}\ \underset{n\rightarrow \infty 
}{\text{liminf}}\sup_{h\in H_{n,b/2,2}^{\ast }}\mathbf{E}_{h}\left[ \tau
_{M}\left( |\sqrt{n}\{\hat{\theta}-f\left( g(\mathbf{\beta }_{n}(h))\right)
|\right) \right]  \label{st2} \\
&\geq &\inf_{c\in \mathbf{R}}\sup_{\mathbf{r}\in \mathbf{R}^{d}}\mathbf{E}%
\left[ \tau \left( \left\vert f_{-}^{\prime }(g(\mathbf{\beta }%
_{0}))\right\vert \left\vert \tilde{g}_{0}\left( Z+\mathbf{r}\right) -\tilde{%
g}_{0}(\mathbf{r})+c\right\vert \right) \right] .  \notag
\end{eqnarray}%
Combining the bounds in (\ref{st1})\ and (\ref{st2}) into (\ref{bd7}), we
conclude that%
\begin{equation}
\underset{b\rightarrow \infty }{\text{lim}}\ \underset{n\rightarrow \infty }{%
\text{liminf}}\sup_{h\in H_{n,b}}\mathbf{E}_{h}\left[ \tau _{M}\left( |\sqrt{%
n}\{\hat{\theta}-f(g(\mathbf{\beta }_{n}(h)))|\right) \right] \geq
\max_{l=1,2}\Psi (a_{l}),  \label{bound}
\end{equation}%
where%
\begin{equation*}
\Psi (a)\equiv \inf_{c\in \mathbf{R}}\sup_{\mathbf{r}\in \mathbf{R}^{d}}%
\mathbf{E}\left[ \tau \left( a\left\vert \tilde{g}_{0}\left( Z+\mathbf{r}%
\right) -\tilde{g}_{0}(\mathbf{r})+c\right\vert \right) \right] ,
\end{equation*}%
and 
\begin{equation*}
a_{+}\equiv \left\vert f_{+}^{\prime }(g(\mathbf{\beta }_{0}))\right\vert 
\text{ and }a_{-}\equiv \left\vert f_{-}^{\prime }(g(\mathbf{\beta }%
_{0}))\right\vert .
\end{equation*}%
Note that $\Psi (a)$ is an increasing function of $a$ on $[0,\infty )$.
Hence the last bound is equal to%
\begin{equation*}
\Psi (\max \left\{ a_{+},a_{-}\right\} ).
\end{equation*}%
Since $f^{\prime }$ is Lipschitz continuous on $[x_{1},g(\mathbf{\beta }%
_{0}))$ and $(g(\mathbf{\beta }_{0}),x_{2}]$ by Assumption 1(ii), we have%
\begin{eqnarray*}
a_{+} &=&\text{lim}_{y\downarrow 0}|f^{\prime }(g(\mathbf{\beta }_{0})+y)|=%
\text{lim}_{\varepsilon \downarrow 0}\sup_{0<y\leq \varepsilon }|f^{\prime
}(g(\mathbf{\beta }_{0})+y)|\text{ and} \\
a_{-} &=&\text{lim}_{y\downarrow 0}|f^{\prime }(g(\mathbf{\beta }_{0})-y)|=%
\text{lim}_{\varepsilon \downarrow 0}\sup_{0<y\leq \varepsilon }|f^{\prime
}(g(\mathbf{\beta }_{0})-y)|.
\end{eqnarray*}%
Since max function is continuous, 
\begin{eqnarray*}
\max \left\{ a_{+},a_{-}\right\} &=&\text{lim}_{\varepsilon \downarrow
0}\max \left\{ \sup_{0<y\leq \varepsilon }|f^{\prime }(g(\mathbf{\beta }%
_{0})+y)|,\sup_{0<y\leq \varepsilon }|f^{\prime }(g(\mathbf{\beta }%
_{0})-y)|\right\} \\
&=&\text{lim}_{\varepsilon \downarrow 0}\sup_{y\in \lbrack -\varepsilon
,\varepsilon ]\backslash \{0\}}|f^{\prime }(g(\mathbf{\beta }_{0})+y)|=\bar{f%
}^{\prime }(g(\mathbf{\beta }_{0})).
\end{eqnarray*}%
Thus we have a desired lower bound. $\blacksquare $\bigskip

For a given $M_{1}>0$, define%
\begin{equation*}
B_{M_{1}}(c)\equiv \sup_{\mathbf{r}\in \mathbf{R}^{d}}\mathbf{E}\left[ \tau
_{M_{1}}\left( a_{0}\left\vert \tilde{g}_{0}(Z+\mathbf{r})-\tilde{g}_{0}(%
\mathbf{r})+c\right\vert \right) \right] ,
\end{equation*}%
where $a_{0}\equiv \bar{f}^{\prime }(g(\mathbf{\beta }_{0}))$, and let%
\begin{equation*}
E_{M_{1}}\equiv \left\{ c\in \lbrack -M_{1},M_{1}]:B_{M_{1}}(c)\leq
\inf_{c_{1}\in \lbrack -M_{1},M_{1}]}B_{M_{1}}(c_{1})\right\} .
\end{equation*}%
Define $c_{M_{1}}\equiv 0.5\left\{ \max E_{M_{1}}+\min E_{M_{1}}\right\} $.
We also define 
\begin{equation*}
\bar{g}_{n}(\mathbf{z})\equiv g\left( \mathbf{z+}\varepsilon _{n}^{-1}(%
\mathbf{\beta }_{0}-g(\mathbf{\beta }_{0}))\right) ,
\end{equation*}%
for $\mathbf{z\in R}^{d}$, and%
\begin{eqnarray*}
\bar{B}_{M_{1}}(c) &\equiv &\sup_{\mathbf{r}\in \lbrack -M_{1},M_{1}]^{d}}%
\frac{1}{L}\sum_{i=1}^{L}\tau _{M_{1}}\left( \hat{a}_{n}\left\vert \bar{g}%
_{n}(\hat{\Sigma}^{1/2}\mathbf{\xi }_{i}+\mathbf{r})-\bar{g}_{n}(\mathbf{r}%
)+c\right\vert \right) \text{,} \\
\tilde{B}_{M_{1}}(c) &\equiv &\sup_{\mathbf{r}\in \lbrack -M_{1},M_{1}]^{d}}%
\frac{1}{L}\sum_{i=1}^{L}\tau _{M_{1}}\left( a_{0}\left\vert \bar{g}%
_{n}(\Sigma ^{1/2}\mathbf{\xi }_{i}+\mathbf{r})-\bar{g}_{n}(\mathbf{r}%
)+c\right\vert \right) \text{,}
\end{eqnarray*}%
and%
\begin{equation*}
B_{M_{1}}^{\ast }(c)\equiv \sup_{\mathbf{r}\in \lbrack -M_{1},M_{1}]^{d}}%
\mathbf{E}\left[ \tau _{M_{1}}\left( a_{0}\left\vert \bar{g}_{n}(\Sigma
^{1/2}\mathbf{\xi }_{i}+\mathbf{r})-\bar{g}_{n}(\mathbf{r})+c\right\vert
\right) \right] .
\end{equation*}%
We also define%
\begin{equation*}
E_{M_{1}}^{\ast }\equiv \left\{ c\in \lbrack -M_{1},M_{1}]:B_{M_{1}}^{\ast
}(c)\leq \inf_{c_{1}\in \lbrack -M_{1},M_{1}]}B_{M_{1}}^{\ast
}(c_{1})\right\} .
\end{equation*}

\noindent \textbf{Lemma A1:} \textit{Suppose that Assumptions} 1(i), 2,%
\textit{\ and }3 \textit{hold.}

\noindent (i) \textit{As }$K\rightarrow \infty ,$%
\begin{equation*}
\lim_{n\rightarrow \infty }\sup_{h\in H}P_{n,h}\left\{ \sup_{c\mathbf{\in }%
[-M_{1},M_{1}]}\left\vert B_{M_{1}}^{\ast }(c)-\hat{B}_{M_{1}}(c)\right\vert
>K(L^{-1/2}+n^{-1/2}\varepsilon _{n}^{-1}+\varepsilon _{n})\right\}
\rightarrow 0.
\end{equation*}%
\noindent (ii) \textit{As }$n\rightarrow \infty ,$%
\begin{equation*}
\sup_{c\mathbf{\in }[-M_{1},M_{1}]}\left\vert B_{M_{1}}^{\ast
}(c)-B_{M_{1}}(c)\right\vert \rightarrow 0.
\end{equation*}

\noindent \textbf{Proof:} (i) As shown in the proof of Lemma A5 of Song
(2014), we find that as $K\rightarrow \infty ,$%
\begin{equation*}
\lim_{n\rightarrow \infty }\sup_{h\in H}P_{n,h}\left\{ \sup_{\mathbf{z}\in 
\mathbf{R}^{d}}\left\vert \bar{g}_{n}(\mathbf{z})-\hat{g}_{n}(\mathbf{z}%
)\right\vert >Kn^{-1/2}\varepsilon _{n}^{-1}\right\} \rightarrow 0.
\end{equation*}%
Therefore, as $K\rightarrow \infty ,$%
\begin{equation*}
\lim_{n\rightarrow \infty }\sup_{h\in H}P_{n,h}\left\{ \sup_{c\in \lbrack
-M_{1},M_{1}]}\left\vert \bar{B}_{M_{1}}(c)-\hat{B}_{M_{1}}(c)\right\vert
>Kn^{-1/2}\varepsilon _{n}^{-1}\right\} \rightarrow 0.
\end{equation*}

Also, for any $\tilde{\varepsilon}_{n}\downarrow 0$ such that $\tilde{%
\varepsilon}_{n}/\varepsilon _{n}\rightarrow 0$ and $\tilde{\varepsilon}_{n}%
\sqrt{n}\rightarrow \infty $ as $n\rightarrow \infty $, we have that 
\begin{equation*}
\inf_{h\in H}P_{n,h}\left\{ \left\vert g(\mathbf{\tilde{\beta})-}g(\mathbf{%
\beta }_{0}\mathbf{)}\right\vert \leq C\tilde{\varepsilon}_{n}\right\}
\rightarrow 1,
\end{equation*}%
for some $C>0$. Therefore, with probability approaching one (uniformly over $%
h\in H$),%
\begin{eqnarray*}
\left\vert \hat{a}_{n}-a_{0}(\tilde{\varepsilon}_{n})\right\vert &\leq
&\sup_{x\in \lbrack g(\mathbf{\tilde{\beta}})-\varepsilon _{n},g(\mathbf{%
\tilde{\beta}})+\varepsilon _{n}]\cap \mathcal{Y}}|f^{\prime }(x)|\  \\
&&-\sup_{x\in \lbrack g(\mathbf{\tilde{\beta}})-\tilde{\varepsilon}_{n},g(%
\mathbf{\tilde{\beta}})+\tilde{\varepsilon}_{n}]\cap \mathcal{Y}}\left\vert
f^{\prime }(x)\right\vert \\
&\leq &\sup_{x\in \lbrack g(\mathbf{\tilde{\beta}})-\varepsilon _{n},g(%
\mathbf{\tilde{\beta}})+\varepsilon _{n}]\cap \mathcal{Y}}|f^{\prime }(x)|\ 
\\
&&-\ \sup_{x\in \lbrack g(\mathbf{\beta }_{0})-\tilde{\varepsilon}_{n}/2,g(%
\mathbf{\beta }_{0})+\tilde{\varepsilon}_{n}/2]\cap \mathcal{Y}}\left\vert
f^{\prime }(x)\right\vert \\
&\leq &C\varepsilon _{n}\rightarrow 0,
\end{eqnarray*}%
as $n\rightarrow \infty $, for some constant $C>0$. The last bound $%
C\varepsilon _{n}$ follows from the assumption that the derivative $%
f^{\prime }(x)$ is Lipschitz continuous on $\mathcal{Y}$. Following the
proof of Lemma A5 of Song (2014), we conclude that%
\begin{eqnarray*}
\lim_{n\rightarrow \infty }\sup_{h\in H}P_{n,h}\left\{ \sup_{c\in \lbrack
-M_{1},M_{1}]}\left\vert \tilde{B}_{M_{1}}(c)-\bar{B}_{M_{1}}(c)\right\vert
>K\{n^{-1/2}+\varepsilon _{n}\}\right\} &\rightarrow &0\text{ and } \\
\lim_{n\rightarrow \infty }P\left\{ \sup_{c\in \lbrack
-M_{1},M_{1}]}\left\vert B_{M_{1}}^{\ast }(c)-\tilde{B}_{M_{1}}(c)\right%
\vert >K(L^{-1/2}+n^{-1/2})\right\} &\rightarrow &0,
\end{eqnarray*}%
as $K\rightarrow \infty $. Combining these results, we obtain the desired
result.

\noindent (ii) The proof is precisely the same as the proof of Lemma A6 of
Song (2014). $\blacksquare $\bigskip

The following lemma deals with the discrepancy between $\hat{c}_{M_{1}}$ and 
$c_{M_{1}}$.\bigskip

\noindent \textbf{Lemma A2:} \textit{Suppose that Assumptions 1(i), 2,\ and
3 hold. Then there exists }$M_{0}$\textit{\ such that for any} $M_{1}>M_{0},$%
\textit{\ and any }$\varepsilon >0,$%
\begin{equation*}
\sup_{h\in H}P_{n,h}\left\{ \left\vert \hat{c}_{M_{1}}-c_{M_{1}}\right\vert
>\varepsilon \right\} \rightarrow 0,
\end{equation*}%
\textit{as} $n,L\rightarrow \infty $ \textit{jointly.}\bigskip

\noindent \textbf{Proof of Lemma A2:} The proof essentially modifies that of
Lemma A7 of Song (2014). From the latter proof, it suffice to show (a)\ and
(b) in the proof of Lemma A7 of Song (2014) for our context. We can derive
these using Lemma A1 precisely in the same way. $\blacksquare $\bigskip

\noindent \textbf{Proof of Theorem 2:} Take $\tilde{\varepsilon}%
_{n}\downarrow 0$ such that $\tilde{\varepsilon}_{n}/\varepsilon
_{n}\rightarrow 0$ but $\tilde{\varepsilon}_{n}/\sqrt{n}\rightarrow \infty $%
. Then, observe that 
\begin{equation*}
\inf_{h\in H}P_{n,h}\left\{ \left\vert g(\mathbf{\tilde{\beta})-}g(\mathbf{%
\beta }_{0}\mathbf{)}+\frac{\hat{c}_{M_{1}}}{\sqrt{n}}\right\vert \leq C%
\tilde{\varepsilon}_{n}\right\} \rightarrow 1,
\end{equation*}%
for some constant $C>0$ that does not depend on $h\in H$ by Assumption 3(i)
and Lipschitz continuity of $g$ and the fact that%
\begin{equation*}
\frac{|\hat{c}_{M_{1}}|}{\sqrt{n}}\leq \frac{M_{1}}{\sqrt{n}}\rightarrow 0,
\end{equation*}%
as $n\rightarrow \infty $ for each fixed $M_{1}>0$. Then, with probability
approaching one,%
\begin{eqnarray}
\left\vert \sqrt{n}\{\hat{\theta}_{mx}-f(g(\mathbf{\beta }%
_{n}(h)))\}\right\vert &\leq &a_{n}\left\vert g(\mathbf{\tilde{\beta}})-g(%
\mathbf{\beta }_{n}(h))+\frac{\hat{c}_{M_{1}}}{\sqrt{n}}\right\vert
\label{bd11} \\
&=&a_{n}\left\vert g(\mathbf{\tilde{\beta}})-g(\mathbf{\beta }_{n}(h))+\frac{%
c_{M_{1}}}{\sqrt{n}}\right\vert +o_{P}(n^{-1/2}),  \notag
\end{eqnarray}%
where 
\begin{equation}
a_{n}\equiv \sup_{x\in \lbrack g(\mathbf{\beta }_{0})-\tilde{\varepsilon}%
_{n}/2,g(\mathbf{\beta }_{0})+\tilde{\varepsilon}_{n}/2]\cap \mathcal{Y}%
}|f^{\prime }(x)|,  \label{a_sup}
\end{equation}%
and $o_{P}(1)$ is uniform over $h\in H$. The last equality follows by
Assumption 3(i), Lipschitz continuity of $f$, and Lemma A1. Note that the
supremum in (\ref{a_sup}) is monotone decreasing in $n$, so that%
\begin{equation}
a_{n}\rightarrow \bar{f}^{\prime }(g(\mathbf{\beta }_{0})),  \label{cv21}
\end{equation}%
as $n\rightarrow \infty $.

Therefore, using (\ref{bd11}) and following precisely the same proof as that
of Theorem 2 in Song (2014), we have%
\begin{eqnarray*}
&&\underset{n\rightarrow \infty }{\text{limsup}}\sup_{h\in H_{n,b}}\mathbf{E}%
_{h}\left[ \tau _{M}(|\sqrt{n}\{\hat{\theta}_{mx}-f(g(\mathbf{\beta }%
_{n}(h)))\}|)\right] \\
&\leq &\ \underset{n\rightarrow \infty }{\text{limsup}}\sup_{\mathbf{r}\in 
\mathbf{R}^{d}}\mathbf{E}\left[ \tau _{M}\left( a_{n}\left\vert \tilde{g}%
_{0}\left( Z+\mathbf{r}\right) -\tilde{g}_{0}\left( \mathbf{r}\right)
+c_{M_{1}}\right\vert \right) \right] .
\end{eqnarray*}%
By (\ref{cv21}), the last term is bounded by%
\begin{eqnarray*}
&&\sup_{\mathbf{r}\in \mathbf{R}^{d}}\mathbf{E}\left[ \tau _{M_{1}}\left( 
\bar{f}^{\prime }(g(\mathbf{\beta }_{0}))\left\vert \tilde{g}_{0}\left( Z+%
\mathbf{r}\right) -\tilde{g}_{0}\left( \mathbf{r}\right)
+c_{M_{1}}\right\vert \right) \right] \\
&=&\inf_{c\in \lbrack -M_{1},M_{1}]}\sup_{\mathbf{r}\in \mathbf{R}^{d}}%
\mathbf{E}\left[ \tau _{M_{1}}\left( \bar{f}^{\prime }(g(\mathbf{\beta }%
_{0}))\left\vert \tilde{g}_{0}\left( Z+\mathbf{r}\right) -\tilde{g}%
_{0}\left( \mathbf{r}\right) +c\right\vert \right) \right] .
\end{eqnarray*}%
The last equality follows by the definition of $c_{M_{1}}$. Finally, we
increase $M_{1}\uparrow \infty $ to obtain the desired result. $\blacksquare 
$\bigskip

\noindent \textbf{Proof of Theorem 3:} First, due to Assumption 3', the
convergences in Lemmas A1 and A2 are uniform over $\alpha _{0}\in \mathcal{A}
$. Take $\tilde{\varepsilon}_{n}\downarrow 0$ such that $\tilde{\varepsilon}%
_{n}/\varepsilon _{n}\rightarrow 0$ but $\tilde{\varepsilon}_{n}/\sqrt{n}%
\rightarrow \infty $. Then, observe that by Assumption 3'(i), Lipschitz
continuity of $g$,%
\begin{equation*}
\inf_{\alpha \in \mathcal{A}_{n}(\alpha _{0};\varepsilon _{n})}\inf_{h\in
H}P_{n,h}\left\{ \left\vert g(\mathbf{\tilde{\beta})-}g(\mathbf{\beta }%
(P_{\alpha })\mathbf{)}+\frac{\hat{c}_{M_{1}}}{\sqrt{n}}\right\vert \leq C%
\tilde{\varepsilon}_{n}\right\} \rightarrow 1,
\end{equation*}%
for some constant $C>0$ that does not depend on $h\in H$ or $\alpha \in 
\mathcal{A}$. Then, we have%
\begin{equation*}
\left\vert \sqrt{n}\{\hat{\theta}_{mx}-f(g(\mathbf{\beta }%
_{n}(h)))\}\right\vert \leq \bar{a}_{n}\left\vert g(\mathbf{\tilde{\beta}}%
)-g(\mathbf{\beta }_{n}(h))+\frac{c_{M_{1}}}{\sqrt{n}}\right\vert
+o_{P}(n^{-1/2}),
\end{equation*}%
where $o_{P}(1)$ is uniform over $h\in H$ and over $\alpha \in \mathcal{A}$,
and%
\begin{equation*}
\bar{a}_{n}\equiv \sup_{\alpha \in \mathcal{A}_{n}(\alpha _{0};\varepsilon
_{n})}\sup_{x\in \lbrack g(\mathbf{\beta }_{0}\mathbf{)}-\tilde{\varepsilon}%
_{n},g(\mathbf{\beta }_{0}\mathbf{)}+\tilde{\varepsilon}_{n}]\cap \mathcal{Y}%
}|f^{\prime }(x)|.
\end{equation*}%
(Recall that $\mathbf{\beta }_{0}=\mathbf{\beta }(P_{\alpha _{0}})$ and
hence it depends on $\alpha _{0}\in \mathcal{A}$.) Similarly as in the proof
of Theorem 2,%
\begin{eqnarray*}
&&\underset{n\rightarrow \infty }{\text{limsup}}\sup_{\alpha \in \mathcal{A}%
_{n}(\alpha _{0};\varepsilon _{n})}\sup_{h\in H_{n,b}}\mathbf{E}_{h}\left[
\tau _{M}(|\sqrt{n}\{\hat{\theta}_{mx}-f(g(\mathbf{\beta }_{n}(h)))\}|)%
\right]  \\
&\leq &\ \underset{n\rightarrow \infty }{\text{limsup}}\sup_{\mathbf{r}\in 
\mathbf{R}^{d}}\mathbf{E}\left[ \tau _{M}\left( \bar{a}_{n}\left\vert \tilde{%
g}_{0}\left( Z+\mathbf{r}\right) -\tilde{g}_{0}\left( \mathbf{r}\right)
+c_{M_{1}}\right\vert \right) \right]  \\
&\leq &\ \sup_{\mathbf{r}\in \mathbf{R}^{d}}\mathbf{E}\left[ \tau _{M}\left( 
\bar{a}_{n_{0}}\left\vert \tilde{g}_{0}\left( Z+\mathbf{r}\right) -\tilde{g}%
_{0}\left( \mathbf{r}\right) +c_{M_{1}}\right\vert \right) \right] ,
\end{eqnarray*}%
for any fixed $n_{0}\geq 1$. The last inequality follows because $%
\varepsilon _{n}\downarrow 0$ and $\tilde{\varepsilon}_{n}\downarrow 0$ as $%
n\rightarrow \infty $, and $\bar{a}_{n}$ and $\mathcal{A}_{n}(\alpha
_{0};\varepsilon _{n})$ are decreasing in $n$. We send $n_{0}\rightarrow
\infty $ and apply the monotone convergence theorem to obtain the last bound
as%
\begin{eqnarray*}
&&\sup_{\mathbf{r}\in \mathbf{R}^{d}}\mathbf{E}\left[ \tau _{M_{1}}\left( 
\bar{f}^{\prime }(\mathbf{\beta }_{0})\left\vert \tilde{g}_{0}\left( Z+%
\mathbf{r}\right) -\tilde{g}_{0}\left( \mathbf{r}\right)
+c_{M_{1}}\right\vert \right) \right]  \\
&=&\inf_{c\in \lbrack -M_{1},M_{1}]}\sup_{\mathbf{r}\in \mathbf{R}^{d}}%
\mathbf{E}\left[ \tau _{M_{1}}\left( \bar{f}^{\prime }(\mathbf{\beta }%
_{0})\left\vert \tilde{g}_{0}\left( Z+\mathbf{r}\right) -\tilde{g}_{0}\left( 
\mathbf{r}\right) +c\right\vert \right) \right] .
\end{eqnarray*}%
Finally, we increase $M_{1}\uparrow \infty $ to obtain the desired result. $%
\blacksquare $

\section{Acknowledgement}

This research was supported by the Social Sciences and Humanities Research
Council of Canada.

\end{document}